\documentclass[journal]{IEEEtran}

\usepackage{multibib}
\usepackage{graphicx}
\usepackage{mathptmx} 
\usepackage{newtxtext}
\usepackage{textcomp}
\usepackage[varg,bigdelims]{newtxmath}

\usepackage{enumitem}

\usepackage{amsmath}
\usepackage{dsfont, url}
\usepackage{amsthm}
\usepackage{mathtools}
\usepackage[usenames, dvipsnames]{xcolor}
\usepackage[%
colorlinks={true},%
citecolor={Blue},%
urlcolor={PineGreen},%
linkcolor={Sepia}%
]{hyperref}
\usepackage{float}
\usepackage{newfloat}
\usepackage[skip=1pt]{caption}

\usepackage{amsmath} 
\usepackage{algorithm}
\usepackage[noend]{algpseudocode}
\usepackage{todonotes}
\usepackage{pgf,tikz}
\usepackage{adjustbox}
\usepackage{verbatim} 
\usepackage{xargs}
\usepackage{cite}
\usepackage{multirow}

\usetikzlibrary{plotmarks}
\usetikzlibrary{chains}
\usepackage{pgfplots}
\pgfplotsset{compat=newest}
\pgfplotsset{plot coordinates/math parser=false}
\newlength\figureheight
\newlength\figurewidth
\usetikzlibrary{shapes,arrows}
\usetikzlibrary{positioning, calc, fit}
\usetikzlibrary{decorations,decorations.pathmorphing,decorations.pathreplacing}
\usetikzlibrary{shapes.misc}
\usetikzlibrary{shadows,trees}

\DeclareMathOperator*{\sbjto}{s.\ t.\ }

\DeclareMathOperator{\rank}{rank}
\DeclareMathOperator{\sat}{sat}
\DeclareMathOperator{\trace}{tr}
\DeclareMathOperator{\proj}{Proj}

\renewcommand{\leq}{\leqslant}

\renewcommand{\geq}{\geqslant}

\newcommand{\R}{\mathds{R}}
\newcommand{\Nz}{\mathds{N}_0}
\newcommand{\N}{\mathds{Z}_{+}}

\newcommand{\bmat}[1]{\begin{bmatrix}#1\end{bmatrix}}
\newcommand{\abs}[1]{\left|#1\right|}
\newcommand{\norm}[1]{\left\|#1\right\|}
\newcommand{\secref}[1]{\S \ref{#1}}

\renewcommand{\transp}{^\top}

\newcommand{\zeros}{\mathbf{0}}

\newcommand{\st}{x}

\newcommand{\control}{u}
\newcommand{\controlset}{\mathds{U}}
\newcommand{\wnoise}{w}

\newcommand{\costps}{c_{\mathrm{s}}}
\newcommand{\costfinal}{c_{\mathrm{f}}}

\newcommand{\authority}{u_{\max}}

\newcommand{\Let}{\coloneqq}
\newcommand{\teL}{\eqqcolon}

\newtheorem{assumption}{Assumption}

\newtheorem{theorem}{Theorem}
\newtheorem{lemma}{Lemma}
\newtheorem{proposition}{Proposition}

\newtheorem{definition}{Definition}

\newtheorem{example}{Experiment}

\newtheorem{pstatement}{Problem Statement}

\allowdisplaybreaks
\begin{document}
	
	\title{Deep Model Predictive Control with Stability Guarantees}

	\author{Prabhat K. Mishra,~\IEEEmembership{Member,~IEEE,}
		Mateus V. Gasparino,
		Andres E. B. Velasquez,
		and~Girish Chowdhary,~\IEEEmembership{Senior Member,~IEEE}
		\thanks{ P. K. Mishra, and G. Chowdhary are with Coordinated Science Laboratory,
			M Gasparino and A. E. B. Velasquez are with Agricultural and Biological Engineering,
			University of Illinois at Urbana Champaign (UIUC), 
			USA.
			\tt\{pmishra,mvalve2,andru89,girishc\}@illinois.edu}
	}

	
	\maketitle
	
	\begin{abstract}
		 
		This paper presents a deep learning based model predictive control algorithm for control affine nonlinear discrete time systems with matched and bounded state dependent uncertainties of unknown structure. Since the structure of uncertainties is not known, a deep learning based adaptive mechanism is utilized to mitigate disturbances. In order to avoid any unwanted behavior during the learning phase, a tube based model predictive controller is employed, which ensures satisfaction of constraints and input-to-state stability of the closed-loop states. In addition, the proposed approach guarantees the convergence of states to origin under certain verifiable conditions. To ensure stability and undesirable learning transients, a dual-timescale adaptation mechanism is proposed, where the weights of the last layer of the neural network are updated  each time instant while the inner layers are trained on a slower timescale using training data collected online and selectively stored in a buffer on the basis of singular value maximization criterion. Our results are validated through numerical experiments on wing-rock dynamics. These results indicate that the proposed deep-MPC architecture is effective in learning to control safety critical systems without suffering instability drawbacks.
	\end{abstract}
	
	\begin{IEEEkeywords}
		safety critical systems, deep learning, model predictive control, adaptive control
	\end{IEEEkeywords}

	\section{Introduction}
	\IEEEPARstart{M}{odeling} errors and environmental uncertainties are unavoidable in practice. Therefore, purely model based controllers tend to exhibit unexpected or unwanted behaviors in the real-world. One key solution to this problem is to employ learning-based methods that utilize powerful learning elements such as \emph{deep neural networks} (DNN). Such methods attempt to learn a good model of underlying nonlinear dynamics while the system is in operation in a manner that does not compromise safety and performance. The model could also accommodate change and adapt as the physical system or the operating conditions change.  
	\par There is no doubt that DNN is one of the most powerful function approximation tools available today \cite{lecun2015deep}. Although deep learning has been shown to be successful in solving the complicated control problems \cite{Silver_Nature15, Levine_JMLR16, bojarski2016},  ensuring stability and constraint satisfaction during the learning transients remains a key open problem, especially for \emph{safety critical systems} such as aircraft, autonomous vehicles, multi-agent systems and control of complex chemical processes; see \cite{survey_LMPC, robust_action_governor, safe_rl_Krause, robust_regression, safe_exploration, safe_exploration_Krause, safe_RL_Tomlin} and references therein. 
	\par To address the above challenge, the available domain knowledge in terms of approximate model is utilizing in  \cite{Tomlin_IROS11, NN_collision_avoidance}, along with the learning elements. We refer readers to an excellent survey on safe reinforcement learning \cite{safeRL_survey} and references therein. One key approach is to augment the learning based controller with \emph{model predictive control} (MPC) and related methods to enable online adaptation while guaranteeing safety through constraint satisfaction. In particular, MPC relies on the accuracy of model whereas learning based techniques improve the model through data. Therefore, their proper pairing can bring useful features of both methods while compensating their drawbacks.
	
	\par The constraint satisfaction capability of MPC is combined with function approximation capability of DNN in several recent works; see \cite{Deep_MPC_Gopaluni, Deep_MPC_Lucia_18, rapid_MPC_Borreli, minimax_Deep_MPC, Allgower_MPC_approximate} and references therein. Some newly developed approaches like \emph{system level synthesis} \cite{SLS_safe}, stochastic learning MPC \cite{learning_MPC_Mesbah}, reachable tube learning \cite{deepreach} and iterative learning MPC for autonomous car racing  \cite{LMPC_Borrelli}, employ MPC along with some learning mechanism. A linear tube based MPC is employed along with a \emph{neural network} (NN) in \cite{NNMPC_Findeisen} in which MPC is augmented by NN to deal with nonlinearities arising due to the linearization around the operating point. In order to satisfy constraints on control, first MPC is designed, then NN is constrained to the remaining action space. Alternatively, a robust action governor, which projects the learning based control to the safe action space, is utilized in \cite{robust_action_governor}. In order to exploit full capabilities of both DNN and MPC, some mechanism is needed, which can maintain a training data set online by rich exploration without losing the safety guarantee and where MPC is aware of learning actions to make optimal decisions.  
	
	\par Our main goal in this article is to address this gap by creating a learning based MPC architectures with performance and safety guarantees. When uncertainties are structured, they can be simply represented in terms of (possibly) high dimensional feature basis functions and the learning mechanism acts on the disturbances. These disturbance rejecting actions taken by the learning mechanism are experienced by the MPC controller as additional disturbances. If the learning mechanism eventually rejects the disturbance then MPC can ensure asymptotic convergence of closed-loop states while satisfying the underlying constraints \cite{MWGC}. In this article, we extend the results of \cite{MWGC} for unstructured uncertainties. 
	\par When uncertainties are unstructured, neuro-adaptive controllers are generally utilized in the framework of \emph{model reference adaptive control} (MRAC) \cite{MRAC_Nguyen}. One such method is \cite{HHN_08} in which the hidden layer weights of DNN are randomly assigned and only the output layer is trained by using the MRAC type weight update law. Although adaptation of only the output layer weights gives some acceptable performance, it cannot learn complex features \cite{Chowdhary_GP_14, GPMRAC_18, DMRAC}. In \cite{DMRAC}, authors propose DNN based adaptive architecture for continuous time linear systems in which the output layer weights are updated in real time while the hidden layer weights are trained on a slower time scale using batch updates. Since the training of DNN is a slow process, the separate treatment of the output layer and the hidden layers makes the controller implementable in real time \cite{Joshi_Virdi}. In above references \cite{DMRAC, Joshi_Virdi} constraint satisfaction is missing, which is precisely the main area of thrust of the present paper.  
	
	\par The key contribution of this paper is a provably safe and stable (in terms of Lyapunov based safety and stability certificate) learning based MPC for a class of nonlinear systems that utilizes DNN as learning elements. A mathematically rigorous treatment of Safety and stability of control systems with deep learning in the loop is missing in the literature. To achieve our goal, we begin with the method proposed in \cite{DMRAC} to separate the training of the hidden and output layer weights for control affine non-linear discrete time systems. The output layer weights are updated according to discrete adaptive weight update law. Our crucial addition then is tube based MPC to ensure that the uncertain state trajectory stay within a tube around a nominal reference trajectory in order to satisfy state and control constraints.  In contrast to \cite{NNMPC_Findeisen, robust_action_governor}, we employ bounded activation functions in the outermost layer and project their weights in a bounded set to ensure boundedness of the learning based controller. Then we assign the remaining control authority to MPC for the sake of optimality under constraints. The result is a learning based deep MPC architecture for a broad class of nonlinear systems that is capable of leveraging the power of deep learning without sacrificing stability and constraint satisfaction. Our results are verified in a challenging wing-rock simulation in which we simulate sudden large changes in the uncertainty, and show that our architecture has favorable adaptation and long-term learning properties. 
	\par Our original contributions are summarized below:
	\begin{enumerate}[label = (\alph*), leftmargin = *, nosep]
		\item A deep learning based MPC is presented in the presence of matched and unstructured state dependent uncertainties. We extend the results of \cite{MWGC} to handle unstructured uncertainties, and of \cite{DMRAC} to incorporate hard constraints on state and control for safety critical systems.  
		\item We rigorously study the effect of DNN on adaptive controller in Lemma \ref{lem:adaptation}, the effect of adaptive controller on MPC in Lemma \ref{lem:stability} and adapt the main result of \cite{Mayne_tube_NLMPC} in Proposition \ref{prop:main_tube}. 
		\item Our approach ensures the convergence of states to origin under certain verifiable conditions (Theorem \ref{th:asymptotic}). This result cannot be obtained from \cite{Mayne_tube_NLMPC}.
	\end{enumerate}	 
To the best of our knowledge, this is the first rigorous treatment of deep MPC, which is suitable for real time implementation on safety critical systems.   
	\subsection*{Organization}  
	We present a problem setup in \secref{s:setup}. The formulation of adaptive controller and the MPC controller are given in \secref{s:adaptive} and \secref{s:MPC}, respectively. Our overall algorithm and its stability is discussed in \secref{s:stability}. We validate our theoretical results with the help of a numerical experiments in \secref{s:experiment} and conclude in \secref{s:epilogue}. Our proofs are given in the appendix in a consolidated manner. 	 
	\subsection*{Notations}
	We let $\R$ denote the set of real numbers, $\Nz$ the set of non-negative integers and $\N$ the set of positive integers. For a given vector $v$ and positive (semi)-definite matrix $M \succeq \zeros$, $\norm{v}_M^2$ is used to denote $v \transp M v$. For a given matrix $A$, the trace, the largest eigenvalue, pseudo-inverse and Frobenius norm are denoted by $\trace(A)$, $\lambda_{\max}(A)$, $A^{\dagger}$ and $\norm{A}_F$, respectively. By notation $\norm{A}$ and $\norm{A}_{\infty}$, we mean the standard $2-$norm and $\infty-$norm, respectively, when $A$ is a vector, and induced $2-$norm and $\infty-$norm, respectively, when $A$ is a matrix. A vector or a matrix with all entries $0$ is represented by $\zeros$ and $I$ is an identity matrix of appropriate dimensions. We let $M^{(i)}$ denote the $i^{\text{th}}$ column of a given matrix $M$.
	
	\section{Problem setup}\label{s:setup}
	\par Let us consider a discrete time dynamical system 
	\begin{equation}\label{e:system}
		\st_{t+1} = f(\st_t) + g(\st_t)\left( \control_t + h(\st_t) \right),
	\end{equation}
	where
	\begin{enumerate}[leftmargin = *, nosep, label=(1-\alph*), widest = b]
		\item \label{e:constraints} $\st_t \in \mathcal{X} \subset \R^d$, $ \control_t \in \controlset \Let \{v \in \R^m \mid \norm{v}_{\infty} \leq \authority \}$,
		\item $\mathcal{X}\subset \R^d$ is a compact set, 
		\item system function $f:\R^d \rightarrow \R^d $, control influence function $g:\R^d \rightarrow \R^{d \times m}$ are given Lipschitz continuous functions, 
		\item \label{as: bounds_hg} $h(\st_t)$ is the state dependent matched uncertainty at time $t$ such that $g(\st_t)h(\st_t) \in \mathds{W} \Let \{v \in \R^d \mid \norm{v} \leq w_{\max} \}$, $h$ is continuous, $\norm{g(\st_t)} \leq \delta_g$ for some $\delta_g> 0$ and $\rank(g(\st_t)) = m$ for every $\st_t \in \R^d$. 
	\end{enumerate}  
	The general problem description of this article is as follows:
	\begin{pstatement}\rm{
			Present a stabilizing control framework for \eqref{e:system}, which respects physical constraints \ref{e:constraints}, optimizes a given performance index, and reduces the effect of unstructured uncertainties by using a trainable DNN.}   
	\end{pstatement}

	\par Our proposed solution is based on constraint satisfaction and cost minimization capabilities of MPC, and universal approximation property of neural networks. We break the applied control $\control_t$ such that
	\begin{equation}\label{e:total_conrol} 
		\control_t = \control_t^a + \control_t^m, 
	\end{equation}
	where $\control_t^a$ and $\control_t^m$ are the neuro-adaptive and MPC components, respectively. The MPC controller employs only the nominal dynamics of \eqref{e:system}, which is given below for easy reference
	\begin{equation}\label{e:nominal}
		\st_{t+1} = f(\st_t) + g(\st_t) \control_t^m  \Let \bar{f}(\st_t, \control_t^m).
	\end{equation}
	Therefore, the dynamics \eqref{e:system} can be written as
	\begin{equation}\label{e:actual_agent}
		\st_{t+1} = \bar{f}(\st_t, \control_t^m) + g(\st_t) \left( \control_t^a + h(\st_t)\right).
	\end{equation}
	
	\par Notice that MPC experiences the term $g(\st_t) \left( \control_t^a + h(\st_t)\right)$ as a disturbance. In a broader sense, the MPC component $\control_t^m$ is responsible for input-to-state stability (ISS) of closed-loop states in the presence of bounded disturbances, and the neuro-adaptive control component $\control_t^a$ acts on $h(\st_t)$. In particular, the job of $\control_t^a$ is to approximate $-h(\st_t)$ and keep the approximation error uniformly bounded with a known bound so that MPC can always experience a bounded disturbance. 
	\par In the following sections \secref{s:adaptive} and \secref{s:MPC}, we design the learning based controller $\control_t^a$ and MPC controller $\control_t^m$, respectively.
	\section{Learning based controller}\label{s:adaptive}
	Any continuous function $h$ on a compact set $\mathcal{X}$ can be approximated by a multi-layer network with number of layers $L \geq 2$ such that   
	\begin{equation}
		h(x) = W_L\transp \psi_L \left[ W_{L-1}\transp \psi_{L-1} \left[ \cdots \left [ \psi_1(x)\right] \right] \right] + \varepsilon^{\ast}(x),
	\end{equation} 
	where $x \in \mathcal{X}$, $\psi_i, W_i$ for $i=1, \ldots , L$, are activation functions and ideal weights, respectively, in the $i^{\text{th}}$ layer. The reconstruction error function $\varepsilon^\ast$ is bounded by a known constant $\bar{\varepsilon}^\ast > 0$ for each $x \in \mathcal{X}$, i.\ e.\ $\norm{\varepsilon^\ast(x)} \leq \bar{\varepsilon}^\ast$. Therefore, we can represent $h(\st_t)$ with the help of a neural network with a desired accuracy. If the neural network is not minimal then the ideal weights may not be unique. However, for the neural-adaptive controller design only the existence of ideal weights is assumed, which is always guaranteed when $h$ is a continuous function on a compact set \cite[\S 7.1]{Lewis_NN_99}. Let us define $\phi^{\ast}(x) \Let \psi_L \left[ W_{L-1}\transp \psi_{L-1} \left[ \cdots \left[ \psi_1(x)\right] \right] \right]$ and $W^{\ast} \Let W_{L} $, then	 
	\begin{equation}\label{e:NN_structure}
		h(\st_t) = W^{\ast \top}\phi^{\ast}(\st_t) + \varepsilon^\ast(\st_t),
	\end{equation} 
	where $W^{\ast} \in \R^{(n_L + 1) \times m}$ denotes the weights of output layer. There are $n_L$ number of neurons in the last hidden layer. The first row of $W^{\ast}$ represents the bias term in the output layer and the first $m$ elements of $\phi^\ast$ are $1$. The ideal inner layer weights defining $\phi^{\ast}(\cdot)$ are neither known nor unique. 
	\par We update the weights of the output layer on main machine in real time at each time instant with the help of a weight update law while keeping the weights of hidden layers fixed. The hidden layers are trained on a parallel secondary machine by using the approach \cite{DMRAC} in which the weights of the output layer are copied from the main machine at the start of the training and remain fixed during the training. Once the training of DNN on a secondary machine is complete, new weights of hidden layers are updated on the main machine and remain fixed until new set of weights are again obtained from the secondary machine. 
	\par Let at time $t_0$ the neural network is initialized with random weights on both machines, and for a given $x$ as input, $\phi_0(x)$ denote the output of the last hidden layer at $t_0$. Let $(t_j)_{j \in \N}$ denote the instants when the weights of hidden layers are updated on main machine after the completion of the $j^{\text{th}}$ training. Let $\phi_j(x)$ be the output of the last hidden layer after the $j^{\text{th}}$ training for a given $x$ as input. We can assume that there exist $\varepsilon_{\phi_j}: \R^d \rightarrow \R^m$ for each $\phi_j$ such that $\phi^{\ast}(\st) = \phi_j(\st) + \varepsilon_{\phi_j}(\st)$ for each $\st \in \mathcal{X}$. Therefore, for $t \in \{t_j, t_j +1,  \ldots , t_{j+1}-1 \} $, \eqref{e:NN_structure} becomes  
	\begin{equation}\label{e:NN_structure_general}
		h(\st_t) = W^{\ast \top}\phi_j(\st_t) + \varepsilon_j(\st_t),
	\end{equation} 
	where $\varepsilon_j(\st_t) = \varepsilon^\ast(\st_t) + W^{\ast \top}\varepsilon_{\phi_j}(\st_t)$ is the overall reconstruction error. Notice that even when the hidden layer weights are randomly assigned as in ELM \cite{ELM_06}, the universal approximation property of the neural network allows us to make the overall reconstruction error $\varepsilon_0(\cdot)$ as small as desired by increasing the width of the network. However, a network with trained hidden layers can capture several useful features, which in turn results in performance improvement \cite{DMRAC}.  
	\par We employ 
	\begin{equation}\label{e:adaptive_control}
	\control_t^a = -K_t\transp \phi_j(\st_t)  
	\end{equation}
	as an adaptive (learning) control at time $t \in \{t_j, t_j +1,  \ldots , t_{j+1}-1 \} $, where $K_t$ is the weight of the output layer, which is trained according to the adaptive weight update law and $\phi_j$ is a feature basis function obtained from the inner layers of DNN after $j^{\text{th}}$ training.
	In the next subsections we provide the relevant details of the training of DNN. 
\subsection{Adaptive learning of $W^\ast$ on the main machine}\label{s:outer layer training}
	 We make the following assumption:  
	 \begin{assumption}\label{as:bounds_uncertainty}
	 	\rm{
	 		There exist $\bar{W}_{i} > 0$ for $i = 1, \ldots , m$, and $\sigma, \bar{\varepsilon} > 0$ such that $\norm{W^{\ast (i)}} \leq \bar{W}_{i}$, for $i = 1, \ldots , m$, and $\norm{\phi_j(x)} \leq \sigma, \norm{\varepsilon_j(x)} \leq \bar{\varepsilon}$ for every $x \in \mathcal{X}$ and $j \in \Nz$.}   
	 \end{assumption}
	 The above assumption is standard in literature \cite{HHN_08, NN_15, Joshi_DMRAC}. A priori knowledge about the bounds on the ideal weights $W^\ast$ of the output layer is useful to avoid parameter drift phenomenon. If the activation functions in the last hidden layer are bounded, i.\ e.\ sigmoidal, tanh, etc., then $\norm{\phi_j(x)}$ will also be bounded for each $j$ and for all $x \in \R^d$. 
	 \par We initialize $K_0$ such that $\norm{K_0^{(i)}} \leq \bar{W}_i$; $i = 1, \ldots , m$. For a given learning rate $0< \theta < 1$ and for $t \in \{t_j, t_j +1,  \ldots , t_{j+1}-1 \} $, we employ the following weight update law:
	\begin{equation}\label{e:update_law}
		\bar{K}_{t+1} = K_t + \frac{\theta}{\norm{\phi_j(\st_t)}^2} \phi_j(\st_t)\left( g(\st_t)^{\dagger}(\st_{t+1} - \bar{f}(\st_t, \control_t^m)) \right)\transp , 
	\end{equation}
	where $g(\st_t)^{\dagger} = \left( g(\st_t)\transp g(\st_t) \right)^{-1}g(\st_t)\transp $ represents pseudo-inverse of the left invertible matrix $g(\st_t)$. Notice that first $m$ number of elements in $\phi_j(\cdot)$ are one. Therefore, $\norm{\phi_j(x)}^2 \geq m^2 \geq 1$ for all $x \in \R^d$ and $j \in \Nz$, which avoids any possibility of division by zero. 
	\par We employ the discrete projection method to ensure boundedness of $K_t^{(i)}$ for $i= 1, \ldots , m$, as follows:
	\begin{equation}\label{e:projection}
		K_t^{(i)} = \proj \bar{K}_t^{(i)} = \begin{cases} \bar{K}_t^{(i)} & \text{ if } \norm{\bar{K}_t^{(i)}} \leq \bar{W}_{i} \\
			\frac{\bar{W}_{i}}{\norm{\bar{K}_t^{(i)}}} \bar{K}_t^{(i)} & \text{ otherwise. } \end{cases}
	\end{equation}
	Let $\tilde{K}_t \Let K_t - W^{\ast}$ and $\tilde{\control}_t \Let \control_t^a + h(\st_t) = - \tilde{K}_t \transp \phi_j(\st_t) + \varepsilon_j(\st_t) $. 
	It is evident that $\norm{K_t}_F^2 = \trace(K_t\transp K_t) = \sum_{i=1}^m \norm{K_t^{(i)}}^2 \leq \sum_{i=1}^m \bar{W}_i^2 \teL \bar{W}$ for all $t$ due to the projection. 
	Therefore, the neuro-adaptive control component $\control_t^a $ is bounded, i.\ e.\
	\begin{equation*}
	\norm{\control_t^a} = \norm{K_t\transp \phi(\st_t)} \leq \norm{K_t}\norm{\phi(\st_t)} \leq \norm{K_t}_F \sigma \\
	 \leq \sqrt{\bar{W}}\sigma \teL \authority^a.
	\end{equation*}
	The apparent disturbance term $g(\st_t) \left( \control_t^a + h(\st_t)\right)$ in \eqref{e:actual_agent} is also bounded, i.\ e.\ 
	\begin{equation}\label{e:overall_disturbance_bound}
	\begin{aligned}
	& \norm{g(\st_t) \left( \control_t^a + h(\st_t)\right)} \leq \norm{g(\st_t){\control}_t^a} + \norm{g(\st_t)h(\st_t)} \\
	& \quad \leq \delta_g\norm{\control_t^a} + w_{\max} \leq \delta_g \authority^a + w_{\max} \teL w_{\max}^{\prime}.
	\end{aligned}
	\end{equation}

	\subsection{Semi-supervised learning of $\phi^\ast$ on a secondary machine}\label{s:inner layer training}
	Let $(\bar{t}_j)_{j \in \N}$ represents time instants when we begin the $j^{\text{th}}$ training of DNN. Let $p_0$ number of data samples are required for the training, which are stored in a buffer of size $p_{\max} > p_0$. We do not have access of the labelled data pairs $(x, \phi^\ast(x))$. Therefore, we follow an approach similar to that of \cite{DMRAC} for the data collection and training.
	\par We fix $\bar{t}_1 \geq p_0$ and for each $t \leq \bar{t}_1$, the labelled pairs $(\st_t,\control_t^a)$ are stored in the buffer. Recall that $\control_t^a = -K_t\transp \phi_0(\st_t)$ for $t \leq \bar{t}_1 < t_1$, where $\phi_0(\cdot)$ is obtained by the random initialization of the weights of hidden layers. At $t=\bar{t}_1$, we randomly sample $p_0$ number of data pairs for the training of DNN. We fix the weights of the output layer to be $-K_{\bar{t}_1}$ and train the network. Notice that the training of DNN does not affect the operation of system because the controlled system still employs $\control_t^a = -K_t\transp \phi_0(\st_t)$ as the adaptive control in which only $K_t$ is updated at each time instant by using the weight update law discussed in \secref{s:adaptive}. At $t = t_1$, we get our first trained network. For $t \in \{t_1, \ldots , t_2-1 \}$, we employ $\control_t^a = -K_t\transp \phi_1(\st_t)$ as an adaptive control. This process of training, exploiting and storing is repeated at each time $t$. 
	\par At $t = p_{\max}-1$, the buffer becomes full. So new data can be added after the removal of some old data. We follow the approach of \cite{svd_cl} for the inclusion and removal of data pairs based on the singular value maximization criterion. In particular, we construct a matrix $TT\transp$, where $T$ consists of $p_{\max}$ number of labels, and compute its singular values. If the replacement of $i^{\text{th}}$ label by new label gives larger singular values than the old one, then the new data pair is added at the $i^{\text{th}}$ position of the replay buffer.    
	\par In the next section, we provide important adjustments in tube-based MPC relevant to the context of the present article.	
	
	\section{Model predictive controller}\label{s:MPC}
	The discussion in this section closely follows \cite{Mayne_tube_NLMPC} with some differences, which are occurred due to the inclusion of the adaptive controller discussed in \secref{s:adaptive}. Tube based MPC is based on ensuring that the closed-loop states stay within a tube around a reference trajectory. The reference trajectory is obtained by solving a reference governor problem \emph{offline} under the tightened constraints for regulation problems. Once a trackable reference trajectory is obtained by spending only a part of the available control authority, a reference tracking problem without state constraints is solved \emph{online} that utilizes full control authority. 
	\par Constraint tightening in the reference governor allows satisfaction of actual constraints by the actual states.  Knowledge of the exact bound on disturbance is needed to tighten the constraints. Although the disturbance in dynamical system \eqref{e:system} at time $t$ is $g(x_t)h(x_t)$, the disturbance experienced by MPC is $g(\st_t)(\control_t^a + h(x_t))$, which is uniformly bounded due to \eqref{e:overall_disturbance_bound}. Therefore, we re-define the disturbance set for the reference governor as follows:
	\[ \mathds{W}^{\prime} \Let \{v \in \R^d \mid \norm{v} \leq w_{\max}^{\prime} \} .\]
	We also notice that a part of control is utilized by the adaptive controller. Therefore, we need to re-define the control set as well.
	\[ \controlset^{\prime} \Let \{v \in \R^m \mid \norm{v}_{\infty} \leq \authority - \authority^a \}. \]
	These modifications in tube-based MPC are already pointed out in \cite{MWGC, NNMPC_Findeisen}, for similar context and in \cite{Bhasin_AMPC_19} for the linear systems with parametric uncertainties. 
	\par For some optimization horizon $N \in \N$, an offline reference governor is utilized to generate a reference trajectory 
	\begin{equation}\label{e:reference_signal}
		\begin{aligned}
			(\st_t^r)_{t \in \Nz} &\Let \{ (\st_t^r)_{t = 0}^{N-1}, 0, \ldots \}, \\
			(\control_t^r)_{t \in \Nz} &\Let \{ (\control_t^r)_{t = 0}^{N-1}, 0, \ldots \} .
		\end{aligned}
	\end{equation}
	In particular, the reference trajectory \eqref{e:reference_signal} is obtained by solving the following optimal control problem with penalty matrices $Q,R \succ 0$ and tightened sets $\mathcal{X}_r, \controlset_r$:
	\begin{equation}\label{e:reference_governor}
	\begin{aligned}
	\min_{(\control_i^r)_{i=0}^{N-1}} & \quad \sum_{i=0}^{N-1} \norm{\st_i^r}^2_Q + \norm{\control_i^r}_R^2  \\
	\sbjto \quad  & \st_0^r = \st_0, \st_N^r = \zeros, \\
	&\st_{i+1}^r = \bar{f}(\st_i^r, \control_i^r), \st_i^r \in \mathcal{X}_r \subset \mathcal{X},  \\
	& \control_i^r \in \controlset_r \subset \controlset^{\prime}; i = 0, \ldots, N-1.  
	\end{aligned}
	\end{equation}
	The tightened constraint sets $\mathcal{X}_r$ and $\controlset_r$ can be obtained by following the approach of \cite[\S 7]{Mayne_tube_NLMPC}. 
	\par In order to design the online reference tracking MPC, we first choose an optimization horizon $N \in \N$ and positive definite matrices $Q, R \succ \zeros$, which can be different from those chosen for the reference governor. Let 
	\[ \costps(\st_{t+i \mid t}, \control_{t+i \mid t}) \Let \norm{\st_{t+i \mid t} - \st_{t+i}^r}_Q^2 + \norm{\control_{t+i \mid t} - \control_{t+i}^r}_R^2\]
	be the cost per stage at time $t+i$ predicted at time $t$ and let $\costfinal(x) \Let x\transp Q_f x$ be the terminal cost with $Q_f \succ 0$. The terminal cost $\costfinal$ is treated as a local control Lyapunov function within a terminal set 
	\begin{equation}\label{e:terminal_set}
	\mathcal{X}_f \Let \{x \in \R^d \mid \costfinal(x) \leq \alpha; \; \alpha > 0 \}
	\end{equation}
	 as in \cite{Mayne_tube_NLMPC} by making the following assumption:
	\begin{assumption}\label{as:stability}
		\rm{
			There exists a control $\control^{\prime} \in \controlset^{\prime} $ such that the following holds
			\begin{equation}
				\costfinal \left( \bar{f}(x, \control^{\prime}) \right) -\costfinal(x)  \leq -\costps(x, \control^{\prime})
			\end{equation}
			for every $x \in \mathcal{X}_f$.
		}
	\end{assumption}
The above assumption is standard in literature. Refer to \cite[\S 4]{Mayne_tube_NLMPC} for more details and a minor modification, which we made for simplicity.
The terminal set $\mathcal{X}_f$ is defined as a level set of terminal cost. We will highlight its implication at the end of this section. Let us define 
\begin{equation}\label{e:cost_function}
V_m(\st_{t\mid t}, (\control_{t+i\mid t})_{i=0}^{N-1}) \Let \costfinal(\st_{t+N \mid t}) + \sum_{i=0}^{N-1}\costps(\st_{t+i \mid t}, \control_{t+i \mid t}) .
\end{equation}  
The online reference tracking MPC minimizes \eqref{e:cost_function} at each time instant $t$ under the following constraints:
\begin{align}
&\st_{t \mid t} = \st_t \label{e:constraint_initial}\\
& \control_{t \mid t} + \control_t^a \in \controlset \label{e:constraint_first_control}\\
& \st_{t+ i +1 \mid t} = \bar{f}(\st_{t+i \mid t}, \control_{t+i\mid t}) \text{ for } i = 0, \cdots, N-1   \label{e:constraint_dynamics}\\
& \control_{t + i\mid t} \in \controlset^\prime \text{ for } i = 1, \cdots, N-1 . \label{e:constraint_remaining_control}
\end{align}
We define the underlying optimal control problem as follows:
\begin{equation}\label{e:MPC}
V_m(\st_{t}) \Let \min_{(\control_{t+i \mid t})_{i=0}^{N-1}} \left\{ \eqref{e:cost_function} \mid \eqref{e:constraint_initial}, \eqref{e:constraint_first_control}, \eqref{e:constraint_dynamics}, \eqref{e:constraint_remaining_control} \right\}.
\end{equation}
	Let the optimizer of the above problem be $(\control_{t+i\mid t}^\ast)_{i=0}^{N-1}$. Then optimal cost will be \[ V_m(\st_t) \Let V_m(\st_t, (\control_{t+i\mid t}^\ast )_{i=0}^{N-1} ). \]
	The first control $\control_{t\mid t}^\ast$ is called the MPC component $\control_t^m$ and is applied along with $\control_t^a$ to the system at time $t$, which in turn results in the satisfaction of control constraints. Notice that the \emph{online reference tracking} problem does not have state constraints and it does not take into account the modified disturbance set $\mathds{W}^{\prime}$.

	\section{Algorithm and stability}\label{s:stability} 
	In this section, we present our algorithm and provide our main result on stability. For a regulation problem, the reference governor is utilized at the beginning to generate a trackable reference trajectory. In case of reference tracking problem, the same reference governor can be invoked whenever the target state changes. At the beginning, the weights of hidden layers are randomly assigned and the weights of output layer are set to zero. A finite size replay buffer is used to store the state and adaptive action pair $(\st_t, \control_t^a)$, where $\control_t^a$ is computed according to \secref{s:outer layer training}. The feature basis function $\phi$ needed to compute $\control_t^a$ is time varying and represented by the hidden layers of DNN. Once the sufficient data is stored in the replay buffer, a random subset of it is used to train the DNN with the help of generative network as described in \secref{s:inner layer training}. The replay buffer is updated by singular value maximization approach \cite{svd_cl} and DNN is trained after each fixed interval of time. Once the training of DNN is complete, $\phi$ is also updated. The adaptive control component $\control_t^a$ is communicated to MPC, which solves \eqref{e:MPC} to compute $\control_t^m$. The control $\control_t = \control_t^a + \control_t^m$ is applied to the system. The same procedure is repeated after measuring the next state.      
	Our algorithm is summarized in Algorithm \ref{a:algo}.  
	\begin{algorithm} 
		\caption{deep MPC}
		\label{a:algo}
		\begin{algorithmic}[1]
			\Require  $\st_0 , (t_j)_{j \in \N}$
			\State choose $\theta < 1$ 
			\State Generate the reference trajectory \eqref{e:reference_signal}
			\State initialize $K_0 = \zeros$, $t=0$ and randomly assign weights of hidden layers of DNN to get $\phi_0$ 
			\State initialize a generative network with a buffer on a secondary machine according to \secref{s:inner layer training}
			\For{each $t$}
			\State if $t = t_j$ then update $\phi = \phi_j$, otherwise move to the next step
			\State \label{al:first_step}compute $\control_t^a = - K_t \transp \phi(\st_t)$ 
			\State solve \eqref{e:MPC}, set $\control_t^m = \control_{t \mid t}^\ast$ 
			\State apply $\control_t = \control_t^m + \control_t^a$ to the system and measure $\st_{t+1}$
			\State compute $K_{t+1}$ by the weight update law \eqref{e:update_law} and \eqref{e:projection}
			\EndFor
		\end{algorithmic}
	\end{algorithm}
	\par We recall the following definition:
	\begin{definition}[\cite{Ioannou_Fidan_tutorial}, page 117]\label{def:small_disturbance}
		\rm{
			The vector sequence $(s_t)_{t \in \Nz}$ is called $\mu$ small in mean square sense if it satisfies $\sum_{t=k}^{k+N-1}\norm{s_t}^2 \leq Nc_0 \mu + c_0^{\prime}$ for all $k \in \N$, a given constant $\mu \geq 0$ and some $N \in \N$, where $c_0, c_0^{\prime} \geq 0$. } 
	\end{definition}
	Some straightforward arguments as in \cite[\S 4.11.3]{Ioannou_Fidan_tutorial} give us the following result:
	
	\begin{lemma}\label{lem:adaptation}
		Consider the dynamical system \eqref{e:system}, weight update law \eqref{e:update_law} and the projection method \eqref{e:projection}. Let the Assumption \ref{as:bounds_uncertainty} hold and define $V_a(K_t) \Let \frac{1}{\theta}\trace(\tilde{K}_t\transp\tilde{K}_t)$. Then for all $t$,
		\begin{enumerate}[label={(\rm \roman*)}, leftmargin=*, widest=3, align=left, start=1]
			\item \label{lem:bound_Va} $V_a(K_t) \leq \frac{4}{\theta} \bar{W}$, 
			\item \label{lem:drift_Va} $V_a(K_{t+1}) - V_a(K_t)  \leq - \frac{1- \theta}{\sigma^2} \norm{\tilde{\control}_t}^2 +  \norm{\varepsilon(\st_t)}^2$,
			\item \label{lem:small_tildeu} $\tilde{\control}_t$ is $\bar{\varepsilon}^2$ small in mean square sense with $c_0 = \frac{\sigma^2}{1-\theta }$ and $c_0^{\prime} = \frac{4c_0}{\theta}\bar{W}$ as per the Definition \ref{def:small_disturbance}.	
		\end{enumerate}			
	\end{lemma}
	We provide a proof of Lemma \ref{lem:adaptation} in the appendix. Let $X_c(\st_t^r)$ be the level set around $\st_t^r$ of radius $c$ generated by $V_m(\st_t)$ and $X_c$ be their union. In particular,	
	\begin{equation}\label{e:tube}
	\begin{aligned}
	X_c(\st_t^r) & \Let \{x_t \in \R^d \mid V_m(x_t) \leq c; \; c > 0  \}, \\
	X_c & \Let \cup_{t=0}^{N} X_c(\st_t^r) .
	\end{aligned}
	\end{equation}
	Properties of the value function are summarized in the following Lemma. These results are standard in literature \cite{ref:rawlings-09}. We provide their proofs in the appendix for completeness. 
	\begin{lemma}\label{lem:value_mpc}
		\begin{enumerate}[label={(\rm \roman*)}, leftmargin=*, widest=3, align=left, start=1]
			\item \label{e:satisfaction_terminal_set} If $\alpha \geq c$ then $\st_{t+N\mid t} \in \mathcal{X}_f$ for every $\st_t \in X_c(\st_t^r) $. 
			\item \label{e:iss_Lyapunov1} \cite[Lemma 3]{MWGC} There exist $c_2 > c_1 > 0$ such that 
			\begin{align*} 
			& c_1 \norm{\st_t -\st_t^r}^2 \leq \costps(\st_t, \control_{t\mid t}^\ast) \leq V_m(\st_t) \leq c_2 \norm{\st_t -\st_t^r}^2 .
			\end{align*}
		\end{enumerate}
	\end{lemma} 
	Lemma \ref{lem:value_mpc}-\ref{e:satisfaction_terminal_set} ensures the satisfaction of terminal constraint on states just by construction. Refer to \cite[Proposition 1]{Stewart_Rawling_11} and \cite[Proposition 1]{Mayne_tube_NLMPC} for minor differences due to \eqref{e:system}, \eqref{e:MPC} and Assumption \ref{as:stability}. 
	\par For the purpose of analysis, we define an intermediate optimization problem by replacing $\st_t$ in \eqref{e:constraint_initial} by $\st_{t\mid t-1}$. In particular,
	\begin{equation}\label{e:intermediate_problem}
	\hat{V}_m(\st_{t\mid t-1}) \Let \min_{(\control_{t+i \mid t})_{i=0}^{N-1}} \left\{ \eqref{e:cost_function} \mid \st_{t \mid t} = \st_{t \mid t-1}, \eqref{e:constraint_first_control}, \eqref{e:constraint_dynamics}, \eqref{e:constraint_remaining_control} \right\}.
	\end{equation} 
	Notice that we keep $\control_t^a = -K_t\transp\phi_j(\st_t)$ fixed in both problems \eqref{e:MPC} and \eqref{e:intermediate_problem}, respectively, and therefore, follow the following convention:
	\begin{equation}\label{e:convention}
	\hat{V}_m(\st_{t\mid t-1}) \leq c \implies \st_{t\mid t-1} \in \mathcal{X}_c(\st_t^r). 
	\end{equation}
	Important results related to tube MPC are summarized in the following Lemma. Refer to \cite[Proposition 2, Proposition 4]{Mayne_tube_NLMPC} for a detailed discussion. We provide their proofs in the appendix to highlight the adjustments and for completeness.
	\begin{lemma}\label{lem:stability}
		If Assumption \ref{as:stability} is satisfied, then for all $t$ for every $\st_t \in X_c(\st_t^r)$ the following hold.
		\begin{enumerate}[label={(\rm \roman*)}, leftmargin=*, widest=3, align=left, start=1]
			\item \label{e:bound_on_predicted_next} $\hat{V}_m(\st_{t+1 \mid t}) - V_m(\st_t) \leq -\costps(\st_t, \control_{t \mid t}^\ast)$, and $\st_{t+1\mid t} \in \mathcal{X}_c$.
			\item \label{e:inclusion_next_state} $\st_{t+1} \in \mathcal{X}_c(\st_{t+1}^r) + \mathds{W}^\prime$.
			\item \label{e:actual_predicted} $V_m(\st_{t+1}) - \hat{V}_m(\st_{t+1 \mid t}) \leq c_3 \norm{g(\st_t)\tilde{\control}_t}$.
			\item \label{e:ISS}There exists $\gamma < 1$ such that
			\begin{equation*}
			V_m(\st_{t+1}) \leq \gamma V_m(\st_t) +  c_3 \norm{g(\st_t)\tilde{\control}_t}.
			\end{equation*}
		\end{enumerate}
	\end{lemma}
	The Lemma \ref{lem:stability}-\ref{e:ISS} along with Lemma \ref{lem:value_mpc}-\ref{e:iss_Lyapunov1} ensures that the controlled system is input-to-state stable (ISS) because it admits $V_m$ as an ISS Lyapunov function \cite[Lemma 3.5]{ISS_discrete}. In case of structured uncertainty $\norm{\tilde{\control}_t} \rightarrow 0$ as $t \rightarrow  \infty$, which implies $\norm{\st_t} \rightarrow 0$ \cite[Theorem 1]{MWGC}. Such results are not available in the presence of unstructured uncertainty. However, the existence of invariant and attractive tubes is possible when $\wnoise_{\max}$ and $\authority^a$ are small. We have the following result:  
	\begin{proposition}\label{prop:main_tube}
		Let us define $\bar{c} \Let \frac{c_2 c_3 }{c_1}(\delta_g\authority^a + \wnoise_{\max})$.
		If $\delta_g\authority^a + \wnoise_{\max}  < \frac{c_1}{c_2 c_3}c$, then for all $t\geq N$, the following hold
		\begin{enumerate}[label={(\rm \roman*)}, leftmargin=*, widest=3, align=left, start=1]
			\item \label{prop:invariant_tube} for every $\st_t \in \mathcal{X}_c(0)\setminus \mathcal{X}_{\bar{c}}(0)$, $V_m(\st_{t+1}) < V_m(\st_t)$,
			\item for every $\st_t \in \mathcal{X}_{\bar{c}}(0)$, $\st_{t+1} \in \mathcal{X}_{\bar{c}}(0)$.
			\item In addition, if $\mathcal{X}_c \subset \mathcal{X}$, then $\st_t \in \mathcal{X}$ for all $t$.
		\end{enumerate}
	\end{proposition}
	The Proposition \ref{prop:main_tube} has similar arguments as in \cite[Proposition 4]{Mayne_tube_NLMPC} and confirms the existence of an invariant tube $\mathcal{X}_c(0)$ and an attractive tube $\mathcal{X}_{\bar{c}}(0) \subset \mathcal{X}_c(0)$. 
	\par Suppose there exists some $\hat{N} \geq N$ and $\hat{c}$ such that $\st_{\hat{N}} \in \mathcal{X}_{\hat{c} }(0) \Let \{ x \mid V_m(x) \leq \hat{c} \leq c \}$. Since $c_3$ is a Lipschitz constant of $V_m$ on a compact set $\mathcal{X}_c + \mathds{W}^\prime \supset \mathcal{X}_c(0) \supset \mathcal{X}_{\hat{c}}(0)$, there exists $\hat{c}_3$, which satisfies Lemma \ref{lem:stability}-\ref{e:ISS}. Similarly, let there exists $\hat{\delta}_g \leq \delta_g$ such that $\norm{g(x)} \leq \hat{\delta}_g$ for every $x \in \mathcal{X}_{\hat{c}}(0)$. Since $\bar{c}$ depends on $c_3$ and $\delta_g$, their reduction will result in shrinkage of the attractive tube $\mathcal{X}_{\bar{c}}(0)$. Moreover, since any level set within $\mathcal{X}_{\bar{c}}(0)$ is invariant due to Proposition \ref{prop:main_tube}-\ref{prop:invariant_tube}, a further shrinkage is possible. However, asymptotic convergence is still not guaranteed. If $\gamma^2 < \frac{1}{2}$, then we can get a stronger result provided a certain condition in terms of $c_3$ and $\delta_g$ is satisfied, and the reconstruction error $\varepsilon$ has small gain type property within the invariant tube. We make the following assumption:
	\begin{assumption}\label{as:notation_small_gain}
		  There exists $\beta > 0$ such that $\norm{\varepsilon_j(\st)} \leq \beta \norm{\st}^2$ for all $\st \in \mathcal{X}_{\hat{c}}(0)$ and $j \in \Nz$. 	 
	\end{assumption}
Generally, the norm bound on the reconstruction error is assumed to be linear in $\norm{x}$\cite{HHN_08}. We assumed it to be quadtratic, otherwise the above assumption is standard in literature. 
We have the following result:
\begin{theorem}\label{th:asymptotic}
	Consider the dynamical system \eqref{e:system} controlled by the Algorithm \ref{a:algo}, and let assumptions \ref{as:bounds_uncertainty}, \ref{as:stability} and \ref{as:notation_small_gain} hold. If $\wnoise_{\max}^\prime < \frac{c_1}{c_2 c_3}c$, $\gamma^2 < \frac{1}{2}$ and $\beta < \frac{c_1 m}{\sqrt{2}\sigma \hat{c}_3 \hat{\delta}_g}\sqrt{(1-2\gamma^2)(1-\theta)}$, 
	then $\norm{\st_t} \rightarrow 0$ as $t \rightarrow \infty$.
\end{theorem}
Notice that the main results of tube-based MPC (Proposition \ref{prop:main_tube}) are valid for small disturbances. The Theorem \ref{th:asymptotic} extends Proposition \ref{prop:main_tube} by guaranteeing convergence of states to origin under the conditions on $\gamma$ and $\beta$. Smaller value of $\gamma$ refers to the faster convergence of the value function of nominal MPC. Generally, reconstruction error is comparatively very small with respect to the disturbance. Therefore, the condition on $\gamma$ and $\beta$ are reasonable, and they can be verified in both theoretical and empirical manner.    	    

\section{Numerical experiment}\label{s:experiment}
\par We consider Wing-rock dynamics to corroborate our result. Recall that wing-rock is a benchmark dynamical system, which has been widely used to evaluate adaptive controllers. Unlike existing adaptive control results \cite{yucelen2012command,chowdhary2010concurrent,chowdhary2013concurrent}, we require that the control and system states remain bounded in a pre-specified constraint set. We assume that the given wing-rock dynamical system suffers from bounded and unstructured uncertainties. 
\par Letting $\delta_t$ denote the roll angle in radian, and $p_t$ denote the roll rate in radian per second, the state of the wing-rock dynamics model is $x_t \Let \bmat{\delta_t & p_t}\transp$ at time $t$. We consider the following discrete time dynamics: 
\begin{equation}\label{eq:disc-sys}
x_{t+1} = Ax_t + B\Big(u_t + h(x_t)\Big),
\end{equation}
where $A=\begin{bmatrix}
1 & 0.05 \\ 0 & 1
\end{bmatrix}$, $B=\begin{bmatrix}
0 \\ 0.05
\end{bmatrix}$,
and $h(\cdot)$ is bounded uncertainty. In order to generate $h$ for the purpose of simulation, we use $ h(x_t) = V_t\transp \varsigma(x_t) + \omega_t$, with
$V_t = v_t  V_0, v_t \in \R$, 
\[ V_0\transp = \bmat{0.8 & 0.2314 & 0.6918 & -0.6245 & 0.0095 & 0.0214}, \] and  $\omega_t \in [-\bar{\omega}, \bar{\omega}]$ is a truncated normal random variable with $\bar{\omega} = 0.1523$. 
The function $\varsigma(\cdot)$ is saturated by a standard saturation function as $\varsigma(x) = \sat(\varsigma^{\prime}(x))$, where $\varsigma^{\prime}(x) = \begin{bmatrix}
1 & \delta & p & \abs{\delta} p & \abs{p} p &\delta^3
\end{bmatrix}\transp $ and $\sat(\cdot)$ is a standard saturation function with the threshold $\frac{\bar{\omega}}{5}$.
The controller is not aware of $\varsigma(\cdot)$ and $\omega$. The admissible state and control sets are given below:
\[ \mathcal{X}  = \left[-\frac{\pi}{6},\frac{\pi}{6} \right] \times \left[ -\frac{\pi}{3},\frac{\pi}{3} \right], \text{ and }
\controlset  = \left[-\frac{\pi}{4},\frac{\pi}{4} \right]. \]
\par We compare our proposed approach (deep MPC computed according to the Algorithm \ref{a:algo}) with two controllers, namely tube MPC and shallow MPC. Tube MPC is designed according to \cite{Mayne_tube_NLMPC}. In order to design shallow MPC, we follow our algorithm \ref{a:algo} but we consider only a single layer neural network with $3$ neurons. Therefore, shallow MPC utilizes our adaptive weight update law \eqref{e:update_law} and \eqref{e:projection} to update their weights but the training of NN is not required. 
\par To design the deep MPC, we use a four layer network with sizes $[2, 5, 5, 3]$ respectively, where the first hidden layer has 5 neurons and the outermost layer has $3$ neurons. The weights of the output layer are updated with our adaptive weight update law \eqref{e:update_law} and \eqref{e:projection}, while the remaining three hidden layers are trained using SGD with momentum constant 0.9 and learning rate 0.01. We use nonlinear activation functions after each of the inner layers, and these functions are respectively $[ReLU, ReLU, tanh]$. We use the MATLAB based software package MPT 3.0 \cite{MPT3} to compute the polytopes needed for the constraint tightening. The simulation parameters are chosen to be same in each experiment with optimization horizon $N=10$, $\theta = 0.9, Q = I$, and $R=0.1$. For the reference governor, we used a larger optimization horizon $N=200$, so that the terminal constraint on state can be satisfied.
\begin{example}\label{ex:fixed_multiplier}
	\rm{
In this exeriment, our control objective is to steer the states of the system from $\st_0 = \bmat{\pi/30 & \pi/12}\transp$ to the origin. We set $v_t = 4$ fixed during this experiment. Our experimental results are demonstrated in Fig.\ \ref{fig:experiment1}. 
\par Since the optimization horizon $N=10$ is small, tube based MPC performance is very poor. However, shallow MPC performs reasonably good under the same experimental data. Roll angle trajectory of shallow MPC depicts some offset, whereas deep MPC converges to a very close vicinity of the origin. There are oscillations in the roll rate trajectories due to the presence of non-vanishing random disturbance $\omega_t$.    
		\begin{figure}
			\centering
			\begin{adjustbox}{width = \columnwidth}
				\includegraphics{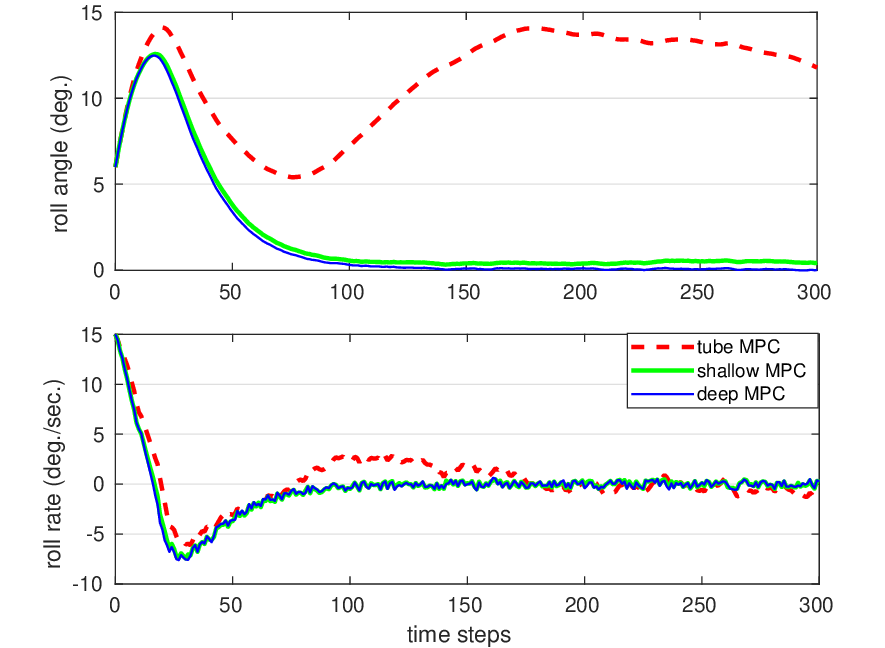}		
			\end{adjustbox}
			\caption{shallow MPC outperforms tube MPC (Experiment \ref{ex:fixed_multiplier})}
			\label{fig:experiment1}
		\end{figure}
	}
\end{example}
\begin{example}\label{ex:regulation}
	\rm{
		In this experiment we consider the system setup same as the Experiment \ref{ex:fixed_multiplier} but the multiplier $v_t$ is considered to be time varying to denote some abrupt changes in uncertainties. The multiplier $v_t$ is given by
		\begin{equation}\label{e:multiplier}
		v_t =  \begin{cases}
		4 \text{ for } t \in [50\ell, 50\ell +49] \quad \ell = 0,2,4, \ldots \\
		0  \text{ otherwise}.
		\end{cases}
		\end{equation}
		\par Our experimental results are depicted in Fig.\  \ref{fig:experiment2}. Due to the sudden change in $v_t$ at time instants shown by vertical grid lines, tube MPC has oscillations in roll angle. Shallow MPC and deep MPC both have similar performance as in the Experiment \ref{ex:fixed_multiplier} with oscillations of comparatively very small magnitude. In this experiment, although the simulation data is same as the Experiment \ref{ex:fixed_multiplier} except the difference in $v_t$, the performance of shallow MPC is affected at each instant of abrupt change (shown by vertical grid line), which depicts its incapability of generalization. However, deep MPC demonstrates a good generalization with only three hidden layers.  
		\begin{figure}
			\centering
			\begin{adjustbox}{width = \columnwidth}
				\includegraphics{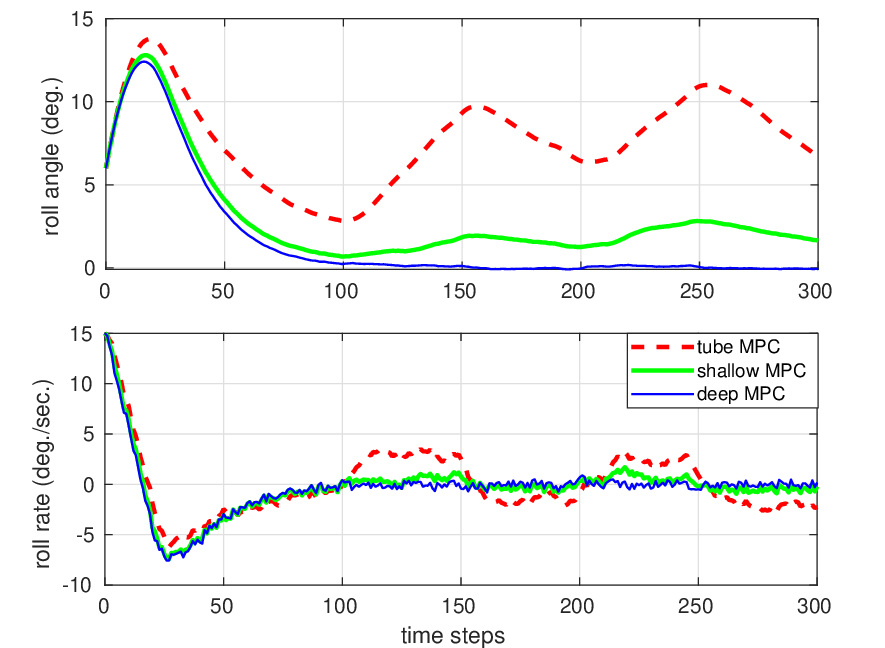}		
			\end{adjustbox}
			\caption{deep MPC outperforms shallow MPC (Experiment \ref{ex:regulation})}
			\label{fig:experiment2}
		\end{figure}
	}
\end{example}
\begin{example}\label{ex:reference-tracking}
	\rm{
		In this experiment, we consider the system \eqref{eq:disc-sys} and multiplier \eqref{e:multiplier} as in the Experiment \ref{ex:regulation}. Our goal is to track a given reference, which is shown in Fig. \ref{fig:tracking} by dotted lines. At each instant, when the reference trajectory changes, we solve the reference governor problem \eqref{e:reference_governor} to generate a trackable reference trajectory.
		\par In order to demonstrate the long term learning of deep MPC, we stopped training of the hidden layers in deep MPC after $t=400$. Fig. \ref{fig:tracking} demonstrates that the proposed deep MPC approach tracks the reference trajectory (shown by dotted lines) very well even after $t=400$, when we fixed the weights of hidden layers in DNN. Since tube MPC and shallow MPC do not have DNN, their behaviors are same as observed in the Experiment \ref{ex:regulation}. 
		\begin{figure*}
			\centering
		\begin{adjustbox}{width = \linewidth}
				\includegraphics{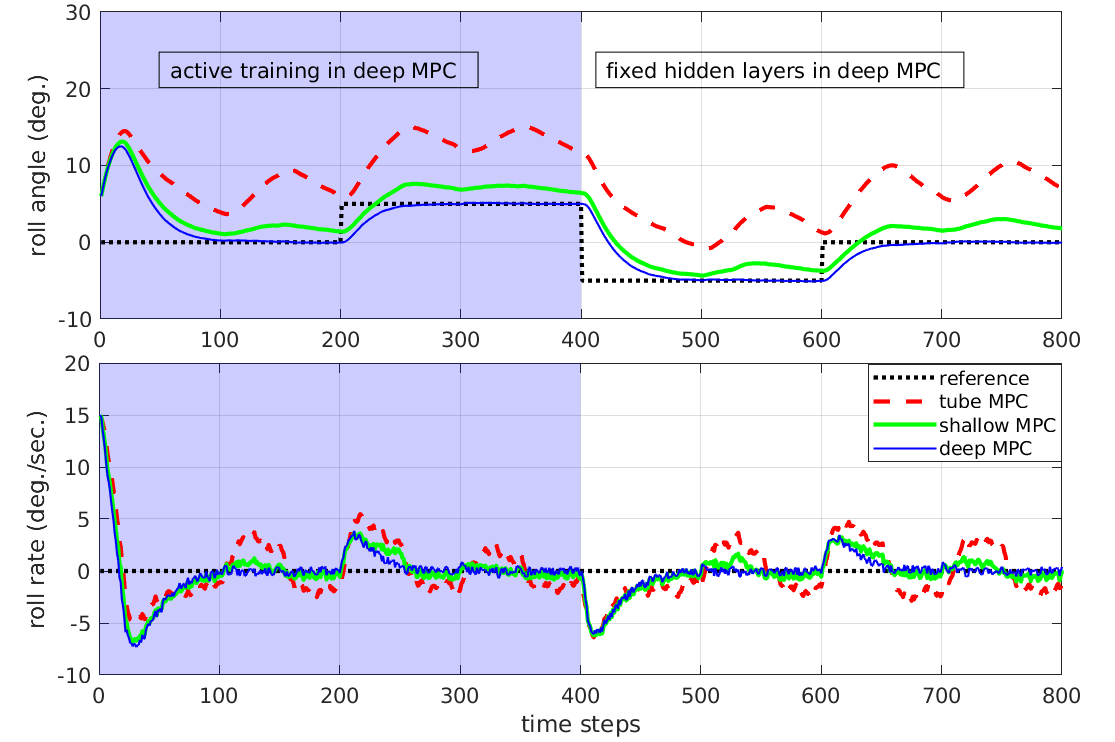}		
			\end{adjustbox}
			\caption{deep MPC tracks a time varying set-point very well even after stopping the training of hidden layers after $400$ time steps (Experiment \ref{ex:reference-tracking})}
			\label{fig:tracking}
		\end{figure*}
	}
\end{example}

\section{Conclusion}\label{s:epilogue}
A deep learning based control algorithm is presented for safety critical systems by combining the approaches of adaptive control, generative neural network and tube MPC. The proposed approach is advantageous over tube MPC when uncertainties are matched. It is demonstrated with the help of numerical experiments that our approach with a single layer neural network (shallow MPC) outperforms tube MPC but results in some offset. The advantage of deep MPC is demonstrated in terms of further improvement in performance and convergence to a very close vicinity of origin. Our approach is conservative in the sense that constraint tightening is required to satisfy constraints. Therefore, it is worth investigating along the lines of \cite{deep_tube_learning, dynamic_tube_How}. Some interesting extensions of the present approach may incorporate output feedback, unreliable channels \cite{PDQ_intermittent, ref:PSC18}. Our approach ensures input-to-state stability of the closed-loop states. In addition, convergence of states to origin is also guaranteed if certain verifiable conditions are satisfied.

		\appendix  

	\begin{proof}[Proof of Lemma \ref{lem:adaptation}]
	\begin{enumerate}[label={(\rm \roman*)}, leftmargin=*, widest=3, align=left, start=1]
	\item Since $V_a(K_t) = \frac{1}{\theta} \trace(\tilde{K}_t\transp \tilde{K}_t) = \frac{1}{\theta} \sum_{i=1}^m \norm{K_t^{(i)} - W^{\ast (i)}}^2 \leq \frac{4}{\theta} \sum_{i=1}^m \bar{W}_i^2 = \frac{4}{\theta} \bar{W}$.
	\item By substituting $\tilde{\control}_t$ in \eqref{e:actual_agent}, we get
	\begin{equation}\label{e:nominal_adaptive}
	\st_{t+1} = \bar{f}(\st_t, \control_t^m) + g(\st_t) \tilde{\control}_t,
	\end{equation}
	and the weight update law \eqref{e:update_law} becomes
	\begin{equation}\label{e:update_law_tilde}
	\begin{aligned}
	\tilde{K}_{t+1} &= \bar{K}_{t+1} - W^\ast + (K_{t+1} - \bar{K}_{t+1}) \\
	&= \tilde{K}_t + \theta \phi_j(\st_t)\tilde{\control}_t \transp + (K_{t+1} - \bar{K}_{t+1}).
	\end{aligned}
	\end{equation}
	Therefore,
	\begin{align*}
	V_a(K_{t+1}) &= \frac{1}{\theta} \trace(\tilde{K}_{t+1}\transp \tilde{K}_{t+1}) \\
	&= \frac{1}{\theta}\trace \left( (\bar{K}_{t+1} - W^\ast)\transp (\bar{K}_{t+1} - W^\ast) + \alpha_t\right),
	\end{align*}
	where 
	\begin{equation*}
	\begin{aligned}
	 \alpha_t &= (K_{t+1} - \bar{K}_{t+1})\transp (K_{t+1} - \bar{K}_{t+1}) \\
	& \quad + 2 (K_{t+1} - \bar{K}_{t+1})\transp (\bar{K}_{t+1} - W^{\ast})\\
    &= -(K_{t+1} - \bar{K}_{t+1})\transp (K_{t+1} - \bar{K}_{t+1})  \\
	& \quad + 2(K_{t+1} - \bar{K}_{t+1})\transp (K_{t+1}-W^\ast). 
	\end{aligned}
	\end{equation*}
	One important property of the projection \eqref{e:projection} is the following \cite[(4.61)]{Ioannou_Fidan_tutorial}:
	\begin{equation}\label{e:effect_projection}
	(W^{\ast (i)} - K_t^{(i)})\transp (\bar{K}_{t}^{(i)} - K_{t}^{(i)})  \leq 0 \text{ for each } i = 1, \ldots , m.
	\end{equation}
	Since $(K_{t+1}^{(i)} - \bar{K}_{t+1}^{(i)})\transp (K_{t+1}^{(i)}-W^\ast) \leq 0$ due to \eqref{e:effect_projection}, we can ensure $\trace(\alpha_t) \leq 0$. Therefore,
	\begin{align*}
	V_a(K_{t+1}) & \leq \frac{1}{\theta}\trace \left( (\bar{K}_{t+1} - W^\ast)\transp (\bar{K}_{t+1} - W^\ast)\right) \\
	&= V_a(K_t) + \frac{1}{\norm{\phi_j(\st_t)}^2} \trace \left( \tilde{\control}_t \tilde{\control}_t\transp + 2 \tilde{K}_t \transp \phi_j(\st_t) \tilde{\control}_t\transp \right).
	\end{align*}
	By substituting $\tilde{K}_t \transp \phi_j(\st_t) = -\tilde{\control}_t + \varepsilon_j(\st_t) $ in the above inequality, we get
	\begin{align*}
	V_a(K_{t+1}) & \leq V_a(K_t) + \frac{1}{\norm{\phi_j(\st_t)}^2} \trace \left( \theta \tilde{\control}_t \tilde{\control}_t\transp + 2 (-\tilde{\control}_t + \varepsilon_j(\st_t)) \tilde{\control}_t\transp \right) \\
	&= V_a(K_t) + \frac{1}{\norm{\phi_j(\st_t)}^2} \left( \theta \norm{\tilde{\control}_t}^2  - 2 \tilde{\control}_t\transp(\tilde{\control}_t - \varepsilon_j(\st_t))  \right) \\
	& = V_a(K_t) + \frac{1}{\norm{\phi_j(\st_t)}^2} \left( (\theta -2) \norm{\tilde{\control}_t}^2  + 2 \tilde{\control}_t\transp \varepsilon_j(\st_t)  \right) \\
	& \leq V_a(K_t) + \frac{1}{\norm{\phi_j(\st_t)}^2} \left( (\theta -1) \norm{\tilde{\control}_t}^2  + \norm{\varepsilon_j(\st_t)}^2  \right) \\
	& = V_a(K_t) - \frac{1- \theta}{\norm{\phi_j(\st_t)}^2} \norm{\tilde{\control}_t}^2 + \frac{1}{\norm{\phi_j(\st_t)}^2} \norm{\varepsilon_j(\st_t)}^2 \\
	& \leq V_a(K_t) - \frac{1- \theta}{\sigma^2} \norm{\tilde{\control}_t}^2 +  \frac{1}{m^2}\norm{\varepsilon_j(\st_t)}^2, \\
	\end{align*}	
	where the last inequality is due to $m^2 \leq \norm{\phi(\st_t)}^2 \leq \sigma^2$. Therefore,
	\begin{equation*}
	V_a(K_{t+1}) - V_a(K_t)  \leq - \frac{1- \theta}{\sigma^2} \norm{\tilde{\control}_t}^2 +  \frac{1}{m^2}\norm{\varepsilon(\st_t)}^2 .
	\end{equation*}
	\item Consider Lemma \ref{lem:adaptation}-\ref{lem:drift_Va} to get
	\begin{align*}
	\frac{1 - \theta }{\sigma^2 }\norm{\tilde{\control}_t}^2 & \leq - V_a(K_{t+1}) + V_a(K_t) + \frac{1}{m^2} \norm{\varepsilon_j(\st_t)}^2 \\
	& \leq - V_a(K_{t+1}) + V_a(K_t) + \frac{\bar{\varepsilon}^2}{m^2}.
	\end{align*}
	By summing from $t=k$ to $k+N-1$ in both sides, we get
	\begin{equation*}
	\begin{aligned}
	&\frac{1 - \theta }{\sigma^2 } \sum_{t=k}^{k+N-1} \norm{\tilde{\control}_t}^2 \leq V_a(K_k) + \frac{N}{m^2} \bar{\varepsilon}^2, \\
	& \quad \leq \frac{4}{\theta} \bar{W} + \frac{N}{m^2} \bar{\varepsilon}^2.  
	\end{aligned}
	\end{equation*}
	Therefore, $\tilde{\control}_t$ is $\bar{\varepsilon}^2$ small in mean square sense with $c_0 = \frac{\sigma^2}{(1-\theta)m^2 }$ and $c_0^{\prime} = \frac{4c_0}{\theta}m^2\bar{W}$ as per the Definition \ref{def:small_disturbance}.	
	\end{enumerate}			 	
	\end{proof}	
\begin{proof}[Proof of Lemma \ref{lem:value_mpc}]
	\begin{enumerate}[label={(\rm \roman*)}, leftmargin=*, widest=3, align=left, start=1]
		\item We recall the definitions of $\mathcal{X}_c(\st_t^r)$ and $\mathcal{X}_f$ from \eqref{e:terminal_set} and \eqref{e:tube}, respectively. Now, it is immediate to notice that $\st_t \in X_c(\st_t^r) \implies V_m(\st_t) \leq c \implies \costfinal(\st_{t+N \mid t}) \leq V_m(\st_t) \leq c \leq \alpha  \implies \st_{t+N \mid t} \in \mathcal{X}_f$.
		\item Since $Q \succ 0$ and $f,g$ are Lipschitz continuous, by \cite[Lemma 3]{MWGC} there exist $c_1, c_2 > 0$ such that
		Lemma \ref{lem:value_mpc}-\ref{e:iss_Lyapunov1} hold. We mention key steps here for completeness. Since $ V_m(\st_t) \geq \costps(\st_t, \control_t^m) \geq \norm{\st_t -\st_t^r}^2_Q$, we can choose $c_1 = \lambda_{\min}(Q)$. 
		\par Let $f, g$ be Lipschitz continuous with Lipschitz constants $L_f$ and $L_g$, respectively. We can notice that \eqref{e:MPC} has no constraints on states and the the constraints on control can be satisfied by $(\control_i^r)_{i=t}^{t+N-1}$ at time $t$. 
\par Let us recall the definition of the cost function \eqref{e:cost_function}, then due to the optimality of $V_m(\st_t)$, we get
		\begin{align*}
		V_m(\st_t) &\leq  V(\st_t, (\control_{t+i}^r)_{i=0}^{N-1}) \\
		& = \st_{t+N\mid t}\transp Q_f \st_{t+N\mid t} + \sum_{i=0}^{N-1} (\st_{t+i \mid t}-\st_{t+i}^r)\transp Q (\st_{t+i \mid t}-\st_{t+i}^r) \\
		& \leq \lambda_{\max}(Q_f)\norm{\st_{t+N\mid t}}^2 + \sum_{i=0}^{N-1} \lambda_{\max}(Q) \norm{\st_{t+i\mid t} - \st_{t+i}^r}^2 .
		\end{align*}
		The above inequality is due to the substitution $ (\control_i)_{i=t}^{t+N-1}= (\control_i^r)_{i=t}^{t+N-1}$. Further,
		\begin{align*}
		\st_{t+i \mid t} - \st_{t+i}^r &= f(\st_{t+i -1 \mid t}) -f(\st_{t+i -1}^r) \\
		& \quad + \left( g(\st_{t+i -1 \mid t})-g(\st_{t+i -1}^r)\right) \control_{t+i-1}^r \\
		\norm{\st_{t+i \mid t} - \st_{t+i}^r} & \leq \left( L_f+L_g\norm{\control_{t+i-1}^r} \right) \norm{\st_{t+i-1 \mid t} - \st_{t+i-1}^r}  \\
		&\leq \bar{L}\norm{\st_{t+i-1 \mid t} - \st_{t+i-1}^r}  \leq \bar{L}^i\norm{\st_{t} - \st_{t}^r},
		\end{align*}
		where $\bar{L} = L_f+L_g\authority^r$. Since $\st_{t+N}^r = \zeros$ for all $t$, there exists $c_2 = \bar{L}^N\lambda_{\max}(Q_f) + \sum_{i=0}^{N-1}\bar{L}^i\lambda_{\max}(Q) > \lambda_{\min}(Q) = c_1$.
	\end{enumerate}
\end{proof}

\begin{proof}[Proof of Lemma \ref{lem:stability}]
	\begin{enumerate}[label={(\rm \roman*)}, leftmargin=*, widest=3, align=left, start=1]
		\item Since $\control_{t+1 \mid t}^\ast \in \controlset^\prime$, we get $\control_{t+1 \mid t}^\ast + \control_{t+1}^a \in \controlset$. Therefore, $\control_{t+1\mid t+1} = \control_{t+1 \mid t}^\ast $ is feasible for \eqref{e:intermediate_problem} at time $t+1$. Since $\control_{t+i+1 \mid t}^\ast \in \controlset^\prime$, for $i=1, \ldots N-2$ the control sequence $\control_{t+i+1 \mid t+1} = \control_{t+i+1 \mid t}^\ast$ is also feasible at time $t+1$ for \eqref{e:intermediate_problem}. Under the above control sequence $\st_{t+N\mid t+1} = \st_{t+N \mid t} \in \mathcal{X}_f$. Therefore, $\control_{t+N \mid t+1} = \control^\prime$ is feasible for some $\control^\prime \in \controlset^\prime$ satisfying the Assumption \ref{as:stability}. In this way, we have constructed a feasible control sequence $(\control_{t+i+1\mid t+1})_{i=0}^{N-1}$ and due to the optimality of $\hat{V}_m(\st_{t+1\mid t})$, we get
		\begin{align*}
		& \hat{V}_m(\st_{t+1 \mid t}) \leq \costfinal(\st_{t+N+1 \mid t+1}) + \sum_{i=0}^{N-1}\costps(\st_{t+1+i \mid t+1}, \control_{t+1+i \mid t+1}) \\
		& = \costfinal(\st_{t+N +1 \mid t+1}) + \costps(\st_{t+N \mid t+1}, \control_{t+N \mid t+1}) \\
		& \quad + \sum_{i=0}^{N-2}\costps(\st_{t+1+i \mid t+1}, \control_{t+1+i \mid t+1}) \\
		& = \costfinal(\st_{t+N + 1 \mid t+1}) + \costps(\st_{t+N \mid t}, \control_{t+N \mid t+1}) \\
		& \quad + \sum_{i=0}^{N-2}\costps(\st_{t+1+i \mid t}, \control_{t+1+i \mid t}) \\
		& = \costfinal(\st_{t+N +1 \mid t+1}) + \costps(\st_{t+N \mid t}, \control^\prime) + \sum_{i=1}^{N-1}\costps(\st_{t+i \mid t}, \control_{t+i \mid t}) \\
		& = \costfinal(\st_{t+N+1 \mid t+1}) + \costps(\st_{t+N \mid t}, \control^\prime) - \costfinal(\st_{t+N \mid t}) \\
		& \quad - \costps(\st_{t}, \control_{t\mid t}^\ast) + V_m(\st_t) \\
		& = \costfinal(\bar{f}(\st_{t+N \mid t}, \control^\prime)) + \costps(\st_{t+N \mid t}, \control^\prime) - \costfinal(\st_{t+N \mid t}) \\
		& \quad - \costps(\st_{t}, \control_{t\mid t}^\ast) + V_m(\st_t) \\
		& \leq  - \costps(\st_{t}, \control_{t\mid t}^\ast) + V_m(\st_t) 
		\end{align*}  
		due to the Assumption \ref{as:stability}. Therefore, $\hat{V}_m(\st_{t+1 \mid t}) \leq V_m(\st_t) \leq c$, which implies $\st_{t+1 \mid t} \in \mathcal{X}_c(\st_{t+1}^r) \subset \mathcal{X}_c$ due to our convention \eqref{e:convention}.
		\item Since $\st_{t+1\mid t} = \bar{f}(\st_t, \control_t^m) \in \mathcal{X}_c(\st_{t+1}^r)$, we get $\st_{t+1} = \st_{t+1\mid t} + g(\st_t)\tilde{\control}_t \in \mathcal{X}_c(\st_{t+1}^r) + \mathds{W}^\prime \subset \mathcal{X}_c + \mathds{W}^\prime$ due to the Lemma \ref{lem:stability}-\ref{e:bound_on_predicted_next}.
		\item  Let $(v_{t+i+1})_{i=0}^{N-1}$ be the minimizer of \eqref{e:intermediate_problem}. Since $(v_{t+i+1})_{i=0}^{N}$ satisfies \eqref{e:constraint_first_control} and \eqref{e:constraint_remaining_control}, it is feasible for \eqref{e:MPC} at $t+1$. Therefore, due to the optimality of $V_m(\st_{t+1})$, we get
		\begin{align*}
		&V_m(\st_{t+1}) \leq V_m(\st_{t+1}, (v_{t+i+1})_{i=0}^{N-1}), 
		\end{align*}
		which in turn implies 
		\begin{align*}
		&V_m(\st_{t+1}) - \hat{V}_m(\st_{t+1\mid t}) \leq V_m(\st_{t+1}, (v_{t+i+1})_{i=0}^{N-1}) \\
		& \quad - V_m(\st_{t+1\mid t}, (v_{t+i+1})_{i=0}^{N-1}) \\
		& \quad \leq \abs{V_m(\st_{t+1}, (v_{t+i+1})_{i=0}^{N-1}) - V_m(\st_{t+1\mid t}, (v_{t+i+1})_{i=0}^{N-1}) }. 
		\end{align*}
		Now we notice that the cost function \eqref{e:cost_function} is Lipschitz continuous in its first argument on the set $\mathcal{X}_c + \mathds{W}^\prime$ while keeping the second argument fixed and $\st_{t+1}, \st_{t+1\mid t} \in \mathcal{X}_c + \mathds{W}^\prime$. Since $v_{t+i+1} \in \controlset$ for $i = 0, \ldots, N-1$, there exists some $c_3 > 0$ such that 
		\begin{align*}
		& \abs{V_m(\st_{t+1}, (v_{t+i + 1})_{i=0}^{N-1}) - V_m(\st_{t+1\mid t}, (v_{t+i+1})_{i=0}^{N-1}) }	\\
		& \quad \leq c_3 \norm{\st_{t+1}- \st_{t+1\mid t}} \\
		& \quad = c_3 \norm{g(\st_t)\tilde{\control}_t}.
		\end{align*}
		Since $t$ was arbitrary, the above result holds for all $t$.
		\item We compute a bound on $V_m(\st_{t+1})- V_m(\st_t) = V_m(\st_{t+1}) - \hat{V}_m(\st_{t+1 \mid t}) + \hat{V}_m(\st_{t+1 \mid t})- V_m(\st_{t}) $. Then by combining the results of Lemma \ref{lem:stability}-\ref{e:bound_on_predicted_next} and Lemma \ref{lem:stability}-\ref{e:actual_predicted}, we get $V_m(\st_{t+1})- V_m(\st_t) \leq -\costps(\st_t, \control_{t \mid t}^\ast) + c_3 \norm{g(\st_t)\tilde{\control}_t}$. Then due to Lemma \ref{lem:value_mpc}-\ref{e:iss_Lyapunov1}, we have  
		\begin{align*}
		V_m(\st_{t+1}) \leq \gamma V_m(\st_t) +  c_3 \norm{g(\st_t)\tilde{\control}_t},
		\end{align*} 
		where $\gamma = 1-\frac{c_1}{c_2}< 1$.		
	\end{enumerate}
\end{proof}

\begin{proof}[Proof of Proposition \ref{prop:main_tube}]
	\begin{enumerate}[label={(\rm \roman*)}, leftmargin=*, widest=3, align=left, start=1]
		\item We can observe that $c_3\norm{g(\st_t) \tilde{\control}_t} \leq c_3\wnoise_{\max}^{\prime} = \frac{c_1}{c_2}\bar{c}$. Therefore, due to Lemma \ref{lem:stability}-\ref{e:ISS}, we get
		\begin{align*}
		V_m(\st_{t+1}) & \leq \gamma V_m(\st_t) +  c_3 \norm{g(\st_t)\tilde{\control}_t} \leq \gamma V_m(\st_t) + c_3 \frac{c_1}{c_2}\bar{c} \\
		& = (1-\frac{c_1}{c_2}) V_m(\st_t) + \frac{c_1}{c_2}\bar{c}.  
		\end{align*}
		Since $c \geq V_m(\st_t) > \bar{c}$ for all $\st_t \in \mathcal{X}_c(0)\setminus\mathcal{X}_{\bar{c}}(0)$, we have $V_m(\st_{t+1}) - V_m(\st_t) \leq \frac{c_1}{c_2} (\bar{c} - V_m(\st_t)) < 0$. 
		\item If $V_m(\st_t) \leq \bar{c}$ then $V_m(\st_{t+1}) \leq \gamma V_m(\st_t) +  c_3 \norm{g(\st_t)\tilde{\control}_t} \leq \gamma \bar{c} + \frac{c_1}{c_2}\bar{c} = \bar{c} \implies \st_{t+1} \in \mathcal{X}_{\bar{c}}(0)$.
		\item For every $\st_t \in \mathcal{X}_c(\st_t^r)\subset \mathcal{X}_c \subset \mathcal{X}$, $\st_{t+1} \in \mathcal{X}_c(\st_{t+1}^r) \subset \mathcal{X}_c \subset \mathcal{X}$ due to Proposition \ref{prop:main_tube}-\ref{prop:invariant_tube}.
	\end{enumerate}
\end{proof}
\begin{proof}[Proof of Theorem \ref{th:asymptotic}]
Let us consider $V(\st_t, K_t) \Let V_m^2(\st_t) + a_0V_a(K_t)$, where $a_0 = \frac{2}{1-\theta} \left( \hat{c}_3 \hat{\delta}_g \sigma \right)^2$. Clearly, $V$ is continuous in $\st_t$ and $K_t$, and satisfies:
\begin{equation*}
\frac{a_0}{\theta} \trace(\tilde{K}_t\transp \tilde{K}_t) + c_1^2 \norm{\st_t}^4 \leq V(\st_t, K_t) \leq \frac{a_0}{\theta} \trace(\tilde{K}_t\transp \tilde{K}_t) + c_2^2 \norm{\st_t}^4.
\end{equation*}
for all $t \geq \hat{N} \geq N$. From Lemma \ref{lem:stability}-\ref{e:ISS} we have
\begin{equation*}
V_m^2(\st_{t+1}) \leq 2\gamma^2V_m^2(\st_{t}) + 2\hat{c}_3^2\norm{g(\st_t)\tilde{\control}_t}^2.
\end{equation*} 
Therefore,
\begin{equation}
V_m^2(\st_{t+1}) - V_m^2(\st_{t}) \leq -(1-2\gamma^2)V_m^2(\st_{t}) + 2\hat{c}_3^2 \hat{\delta}_g^2 \norm{\tilde{\control}_t}^2.
\end{equation}
Now, we compute $V(\st_{t+1}, K_{t+1}) - V(\st_t, K_t)$ and substitute $V_a(K_{t+1}) - V_a(K_t) \leq -\frac{1-\theta}{\sigma^2} \norm{\tilde{\control}_t}^2 + \frac{1}{m^2}\norm{\varepsilon_j(\st_t)}^2$ from Lemma \ref{lem:adaptation}-\ref{lem:drift_Va} to get
\begin{align*}
&V(\st_{t+1}, K_{t+1}) - V(\st_t, K_t) \leq -(1-2\gamma^2)V_m^2(\st_{t}) + \frac{a_0}{m^2} \norm{\varepsilon_j(\st_t)}^2 \\
& \quad \quad + 2\hat{c}_3^2 \hat{\delta}_g^2 \norm{\tilde{\control}_t}^2 - a_0 \left(\frac{1-\theta}{\sigma^2} \right) \norm{\tilde{\control}_t}^2 \\
& \quad = -(1-2\gamma^2)V_m^2(\st_{t}) + \frac{a_0}{m^2} \norm{\varepsilon_j(\st_t)}^2 \\
& \quad \leq  -(1-2\gamma^2)c_1^2 \norm{\st_t}^4 + a_0 \left(\frac{\beta}{m} \right)^2 \norm{\st_t}^4 \\
& \quad = -\eta \norm{\st_t}^4,
\end{align*}
where $\eta = (1-2\gamma^2)c_1^2 - \frac{2}{1-\theta } \left( \sigma \hat{c}_3 \hat{\delta}_g\frac{\beta}{m} \right)^2 > 0$ because $\beta < \frac{c_1 m}{\sqrt{2}\sigma \hat{c}_3 \hat{\delta}_g}\sqrt{(1-2\gamma^2)(1-\theta)}$.
Therefore,
\begin{align*}
\norm{\st_t}^4 &\leq \frac{1}{\eta} \left( - V(\st_{t+1}, K_{t+1}) + V(\st_t, K_t) \right)
\end{align*}
By summing from $t=\hat{N}$ to $k+\hat{N}$ on both sides, we get
\begin{align*}
\sum_{t=\hat{N}}^{\hat{N} + k} \norm{\st_t}^4 & \leq  \frac{1}{\eta} \left( - V(\st_{\hat{N}+k+1}, K_{\hat{N}+k+1}) + V(\st_{\hat{N}}, K_{\hat{N}}) \right) \\
& \leq \frac{1}{\eta} V(\st_{\hat{N}}, K_{\hat{N}}) \\
& = \frac{1}{\eta} \left( V_m^2(\st_{\hat{N}}) + a_0 V_a(K_{\hat{N}})\right)\\
& \leq \frac{1}{\eta} \left( \hat{c}^2 +  \frac{4a_0}{\theta} \bar{W} \right), 
\end{align*}
where the last inequality is due to Lemma \ref{lem:adaptation}-\ref{lem:bound_Va} and the fact that $\st_{\hat{N}} \in \mathcal{X}_{\hat{c}}(0)$.
Since the right hand side of the above inequality is independent of $k$, we have $\sum_{t=\hat{N}}^\infty \norm{\st_t}^4 \leq \frac{1}{\eta} \left( \hat{c}^2 +  \frac{4a_0}{\theta} \bar{W} \right)$, which implies $\norm{\st_t} \rightarrow 0$ as $t \rightarrow \infty$.
\end{proof}

	\section*{Acknowledgment}
	We gratefully acknowledge financial support from ONR MURI N00014-19-1-2373 and joint NSF CPS USDA grant 2018-67007-28379.

	\bibliographystyle{IEEEtran}
	\bibliography{IEEEabrv,MPC_Learning}
\end{document}